\documentclass{article} 
\usepackage{amssymb}
\usepackage{amsmath}
\input{epsf.sty}  
\usepackage{color}

\newtheorem{theorem}{Theorem}[section]

\newtheorem{example}{Example}[section]

\newtheorem{proposition}{Proposition}[section]
\newcommand{\ds}{\displaystyle}

\newcommand{\qed}{\hfill\rule{2mm}{2mm}}      
\title{\huge\bf 
 A review on symmetry properties\\ 
 of birth-death processes
}
\author{
{\rm Antonio  Di Crescenzo}\footnote{
Dipartimento di Matematica, Universit\`a degli Studi di Salerno, 
Via Giovanni Paolo II, 132; 84084 Fisciano (SA), Italy. 
email: {\tt adicrescenzo@unisa.it}}, 
\quad 
{\rm Barbara Martinucci}\footnote{
Dipartimento di Matematica, Universit\`a degli Studi di Salerno, 
Via Giovanni Paolo II, 132; 84084 Fisciano (SA), Italy. 
email: {\tt bmartinucci@unisa.it}}
}
\date{\today}
\begin{document}
\maketitle
%
\begin{abstract}\small
In this paper we review some results on time-homogeneous birth-death processes. Specifically, 
for truncated birth-death processes with two absorbing or two reflecting endpoints, we recall the 
necessary and sufficient conditions on the transition rates such that the transition probabilities 
satisfy a spatial symmetry relation. The latter leads to simple expressions for first-passage-time 
densities and avoiding transition probabilities. This approach is thus thoroughly extended to 
the case of bilateral birth-death processes, even in the presence of catastrophes, and to the 
case of a two-dimensional birth-death process with constant rates. 
\\
{\em Keywords:} Truncated processes; Bilateral processes; Transition probabilities; 
Spatial symmetry; Absorption; Reflection; First-passage time; 
Avoiding transition probabilities; Catastrophes.
\\
{\em Mathematical Subject Classification:}   60J80 
\end{abstract}
%
\section{\bf   Introduction}
\label{section:1}
Symmetry is an useful property that is largely employed to study stochastic processes. Many 
investigations deal with invariance to space or time transformations. A classical example is the 
method of images, that has been successfully exploited to obtain results on hitting densities 
of diffusion processes (see, for instance, Daniels \cite{Da82}). Other cases referring to 
one-dimensional and two-dimensional time-homogeneous diffusion processes can be found, 
e.g.\ in Giorno {\em et al.}\ \cite{GiNoRi} and Di Crescenzo {\em et al.}\ \cite{DiCrGiNoRi}, 
respectively. The results provided in such papers are mainly based on the Markov property 
and on certain spatial symmetries of diffusion processes, and on the continuity of their sample 
paths. Hence, the symmetry-based approach can be extended by far to other classes of 
continuous-time Markov processes, such as the birth-death (BD) processes. We recall that 
the continuity of their sample paths on integers is also identified as the `skip-free' property.  
\par
It is well known that BD processes are widely considered within applications in ecology, genetics 
and evolution (cf.\ Crawford and Suchard \cite{CrSu12}, Novozhilov {\em et al.}\ \cite{NoKaKo2005}, 
\cite{NoKaKo2006}, Ricciardi \cite{Ri}), theoretical neurobiology (see Giorno {\em et al.}\ 
\cite{GiLaNoRi}, L\'ansk\'y and Rospars \cite{LaRo93}), chemical physics (see, e.g.\ 
Conolly {\em et al.}\ \cite{CoPaDh}, Flegg {\em et al.}\ \cite{FLPoGr2008}, Keller and 
Valleriani \cite{KeVa2012}), mathematical finance (see Kou and Kou \cite{KoKo2003}), queueing 
(cf.\ Di Crescenzo {\em et al.}\ \cite{DiCrGiKKNo12}, Giorno {\em et al.}\ \cite{GiNeNo}, Lenin 
and Parthasarathy \cite{LePa0}, Parthasarathy and Lenin \cite{PaLe2004}, for instance). 
\par
Due to its relevance in such applications, the first-passage-time (FPT) problem for BD 
processes has been largely investigated. Classical papers devoted to the disclosure of its 
structural properties are Karlin and Mc~Gregor~\cite{KaMc1},~\cite{KaMc2}, 
Keilson~\cite{Ke1},~\cite{Ke2},~\cite{Ke3}, Kijima~\cite{Kij88}, R\"{o}sler~\cite{Ro}, 
Sumita and Masuda~\cite{SuMa}. Other methods in this context are 
based on more complex tools such as combinatorial arguments and spectral analysis. 
The symmetry-based approach has also been successfully exploited in this context. 
\par
On the ground of the above remarks, in this paper we aim to give a short review of the main 
results for symmetric BD processes, with attention to FPT problems and related topics.  
We first recall the conditions on the transition rates leading to a quasi-symmetry property of the 
transition probabilities, and also to simple expressions both for FPT densities and 
for certain relevant avoiding transition probabilities. This is accomplished in Section \ref{section:2} 
for one-dimensional truncated BD processes both with reflecting and absorbing endpoints. 
An indication on the extension to the cases of bilateral processes is also provided. Various 
examples are also discussed in detail.
\par
In Section \ref{section:3} we illustrate the symmetry properties of bilateral BD processes with 
catastrophes. We refer to total catastrophes whose effect is an instantaneous jump of the 
process to state 0. In this case all sample paths going from a negative state to a positive state 
(and vice versa) are forced to pass through 0. This allows to make use of the symmetry 
property and thus to obtain some expressions for FPT densities through state 0 and 
the corresponding avoiding transition probabilities. 
\par
As for diffusion processes, the symmetry-based approach for BD processes can be extended also 
to higher dimensions. In fact, in Section \ref{section:2d} we consider the symmetry properties of 
a two-dimensional  BD process characterized by constant rates. In this case we deal with a spatial 
symmetry in the plane with respect to the straight line $x_2 = x_1 + r$. As for the one-dimensional 
case,  we also discuss the FPT problem and deal with certain avoiding transition probabilities. 
\section{\bf Symmetry properties of truncated processes}\label{section:2}
Let  $\{X(t),t\geq 0\}$ be a one-dimensional BD process and let  ${\cal S}:=\{0,1,\ldots,N\}$ be its state space, 
for some natural number $N> 1$.  According to the terminology adopted in \cite{LePa0} and 
\cite{vaDo1} we refer to $X(t)$ as a truncated BD process. In this section we assume that the 
endpoints $0$ and $N$ are absorbing states, and that $\{1,2,\ldots,N-1\}$ is a communicating class. 
As usual, we denote by $\lambda_n$ and $\mu_n$ the birth and death rates of $X(t)$, 
with $\lambda_0=\mu_0=\lambda_N=\mu_N=0$ and $\lambda_n,\mu_n>0$ 
for $n=1,2,\ldots,N-1$. For all $k,n\in{\cal S}$ the transition probabilities 
\begin{equation}
	p_{k,n}(t)=P\{X(t)=n\,|\,X(0)=k\} 
	\label{equation:pknt} 
\end{equation}
satisfy the following forward equations:
\begin{eqnarray*}
	&& \hspace{-1.6cm} \frac{d}{dt}\,p_{k,0}(t)=\mu_1\,p_{k,1}(t) 	
 \nonumber \\
	&& \hspace{-1.6cm} \frac{d}{dt}\,p_{k,n}(t)=\lambda_{n-1}\,p_{k,n-1}(t)-
	(\lambda_n+\mu_n)\,p_{k,n}(t)+\mu_{n+1}\,p_{k,n+1}(t)
	\qquad (n=1,2,\ldots,N-1) 
 \\
	&& \hspace{-1.6cm} \frac{d}{dt}\,p_{k,N}(t)=\lambda_{N-1}\,p_{k,N-1}(t),	
	\nonumber
\end{eqnarray*}
with initial condition
\begin{equation*}
	\lim_{t\downarrow 0}p_{k,n}(t)=\delta_{k,n},
\end{equation*}
where $\delta_{k,n}$ is the Kronecker's delta. Clearly,  for all $k\in\{1,2,\ldots,N-1\}$, 
$p_{k,0}(t)$ [$p_{k,N}(t)$] is the probability that the BD process has been absorbed 
at $0$ [$N$] up to time $t$, whereas $p_{k,n}(t)$, for  $n\in\{1,2,\ldots,N-1\}$, is the 
probability that the BD process goes from $k$ to $n$ without absorption up to time $t$. 
\par
A spatial symmetry with respect to $N/2$, the mid point of $\cal S$, has been exploited 
in \cite{DiCr98} for $X(t)$. Specifically, hereafter we recall the necessary and sufficient 
conditions on the birth and death  rates $\lambda_n,\mu_n$ such that the ratio of probabilities 
of symmetric sample paths is time-independent (see Theorem 2.1 of  \cite{DiCr98}). 
\begin{theorem}\label{teorsimmassorb} 
Let us set  
\begin{eqnarray*}
	&& x_0=1, 	\nonumber 		\\
	&& x_n=\ds\frac{\mu_1\,\mu_2\cdots\mu_n} 
	{\lambda_{N-1}\,\lambda_{N-2}\cdots\lambda_{N-n}}
	=\ds\frac{\mu_n}{\lambda_{N-n}}\,x_{n-1}
	\qquad (n=1,2,\ldots,N-1),			
	\\
	&& x_N=\ds\frac{\mu_1}{\lambda_{N-1}}\,x_{N-1}. 
	\nonumber
\end{eqnarray*}
Then, for $k=1,2,\ldots,N-1$, the transition probabilities (\ref{equation:pknt}) satisfy 
the quasi-symmetry relation 
\begin{equation}
	p_{N-k,N-n}(t)=\frac{x_n}{ x_k}\,p_{k,n}(t)
	\qquad (n\in {\cal S}; \; t\geq 0)			
	\label{equation:5}
\end{equation}
if and only if  
\begin{eqnarray*}
	&& \lambda_n\,\mu_{n+1}=\lambda_{N-n-1}\,\mu_{N-n} 	
	\qquad (n=1,2,\ldots,N-2),					
 \\
	&& \lambda_n+\mu_n=\lambda_{N-n}+\mu_{N-n}			
	\qquad (n=1,2,\ldots,N-1).			
\end{eqnarray*}
\end{theorem}
\par
The BD process $X(t)$ is said {\em symmetric\/} if its transition probabilities satisfy relation (\ref{equation:5}), 
according to the concept of symmetric continuous-time Markov chain given in Karlin \cite{Ka}. 
Clearly, Eq.\ (\ref{equation:5}) can be also viewed as an extension of the reflection principle for random 
walks (see, e.g.\  \cite{Fe}). Moreover, we observe that if $k$ and $n$ are symmetric states (i.e., 
$k+n=N$) then property (\ref{equation:5}) identifies with the time-reversibility relation (see, 
e.g.\ Karlin and Mc~Gregor~\cite{KaMc2}) 
$$
	p_{n,k}(t)=\frac{\mu_{k+1}\,\mu_{k+2}\cdots\mu_n}{
	\lambda_k\,\lambda_{k+1}\cdots\lambda_{n-1}}\,p_{k,n}(t)
	\qquad (k<n; \; k+n=N). 
$$
\par
The symmetry property given in Theorem \ref{teorsimmassorb} has been successfully 
employed in various FPT problems. Let  
$$
 T^+_{k,s}=\inf\{t>0: X(t)=s\}, \qquad X(0)=k<s\leq N
$$ 
be the upward FPT of $X(t)$ from state $k>0$ to state $s$, and let $g^+_{k,s}(t)$   
be the corresponding probability density function. Since $X(t)$ is a skip-free process, 
i.e.\ only unity jumps are allowed (see \cite{AbWh89}), and due to the Markov property, 
the following renewal equation holds, for $t\geq 0$: 
$$
p_{k,n}(t)=\int_0^t g^+_{k,s}(\vartheta)\,p_{s,n}(t-\vartheta)\,d\vartheta  	
	\qquad (0<k<s\leq n\leq N). 
$$
This allows us to express the upward FPT density for $t\geq 0$ as follows (see Proposition 2.2 
of \cite{DiCr98}): 
\begin{equation}
 g^+_{k,s}(t)=\lambda_{s-1}\Big[p_{k,s-1}(t)-\int_0^t
 g^+_{k,s}(\vartheta)\,p_{s,s-1}(t-\vartheta)\,d\vartheta\Big]
 \qquad (0<k<s<N).
 \label{eq:relfptdensity}
\end{equation}
Roughly speaking,  equation (\ref{eq:relfptdensity}) expresses the FPT density as the difference of the 
probability of all sample paths that exhibit an upward jump from $s-1$ to $s$ close to time $t$ minus the 
probability of the sample paths (among the latter ones) that already passed through state $s$ up to time $t$. 
Similarly as above, denoting by $T^-_{k,s}$  the downward FPT of $X(t)$, and by $g^-_{k,s}(t)$ its 
probability density function, for $0\leq s<k<N$, one has the following relations, for $t\geq 0$: 
$$
 p_{k,n}(t)=\int_0^t g^-_{k,s}(\vartheta)\,p_{s,n}(t-\vartheta)\,d\vartheta  
 \qquad (0\leq n\leq s<k<N), 
$$
$$
g^-_{k,s}(t)=\mu_{s+1}\Big[p_{k,s+1}(t)-\int_0^t
	g^-_{k,s}(\vartheta)\,p_{s,s+1}(t-\vartheta)\,d\vartheta\Big]
	\qquad (0<s<k<N). 
$$
When $N$ is even, the symmetry property exploited in Theorem \ref{teorsimmassorb} allows us to 
express the FPT densities of a symmetric BD process through the mid  point of the state space, 
$s=N/2$, called {\em symmetry state} (see Theorem 2.3 of \cite{DiCr98}).  
\begin{proposition}\label{theoremg}
If $N=2s$, with $s$ integer, then the FPT densities through the symmetry state $s$ of a symmetric BD 
process for $t\geq 0$ are given by  
\begin{eqnarray}
	&& g^+_{k,s}(t) = \lambda_{s-1}\,p_{k,s-1}(t)-\mu_{s+1}\,p_{k,s+1}(t) 
	\qquad (0<k<s),		
	\label{equation:10a}
	\\
	&& g^-_{k,s}(t) = \mu_{s+1}\,p_{k,s+1}(t)-\lambda_{s-1}\,p_{k,s-1}(t) 
	\qquad (s<k<2s).
	\label{equation:10b}
\end{eqnarray}
\end{proposition}
\par  
We remark that from Eq.\ (\ref{equation:10a}) we have 
\begin{equation}
 g^+_{k,s}(t)=\left\{
 \begin{array}{ll}
 o(t), & k=1,2,\ldots, s-3\\
 \lambda_{s-1}\,\lambda_{s-2}\,t +o(t), & k=s-2, \\ 
 \lambda_{s-1}-\lambda_{s-1}(\lambda_{s-1}+\mu_{s-1})\,t +o(t), & k=s-1.
 \end{array}
 \right.
 \label{eq:gsmallt}
\end{equation}
A similar result for $g^-_{k,s}(t)$ can be obtained from (\ref{equation:10b}), or taking into account that 
for a symmetric BD process for $t\geq 0$  one has 
$$
 g^+_{N-k,N-s}(t)=\frac{x_s}{x_k}\, g^-_{k,s}(t) \qquad (0<s<k<N).
$$
\par
For $0<r<N$, with $r$ integer, let ${\cal S}^-_r=\{1,2,\ldots,r-1\}$ and ${\cal S}^+_r=\{r+1,r+2,\ldots,N-1\}$. 
For $n,k\in {\cal S}^-$ or  $n,k\in {\cal S}^+$, let us now introduce the 
{\em avoiding transition probabilities\/} of $X(t)$: 
\begin{equation}
	p_{k,n}^{\langle r \rangle}(t)
	=P\left\{X(t)=n, \; X(\vartheta)\neq r
	\;\hbox{for all}\;\vartheta\in(0,t)\,|\,X(0)=k\right\}.
 \label{eq:pknrt}
\end{equation}
Roughly speaking, Eq.\ (\ref{eq:pknrt}) defines the joint probability that the BD process is in state $n$ 
at time $t$ and that no visit to state $r$ occurred up to time $t$, conditional on $X(0)=k$. We remark 
that the following relation holds: 
$$
 p_{k,n}^{\langle r\rangle}(t)=p_{k,n}(t)-\int_0^t g^{\pm}_{k,r}(\vartheta)\,
	p_{r,n}(t-\vartheta)\,d\vartheta,
	\qquad n, k\in {\cal S}^{\mp}_r.
$$
When $X(t)$ is symmetric we are able to express $p_{k,n}^{\langle r \rangle}(t)$  
in closed form (see Theorem 2.5 of \cite{DiCr98}). 
\begin{proposition}\label{theoremtaboo}
If $N=2s$, with $s$ integer, then the $s$-avoiding transition probabilities of a symmetric BD 
process are given by 
$$
 p_{k,n}^{\langle s \rangle}(t)
 = p_{k,n}(t) - \frac{x_k}{ x_s}\,p_{2s-k,n}(t)  
 \qquad n, k\in {\cal S}^{\pm}_s. 
$$
\end{proposition}
\par
We point out that Propositions \ref{theoremg} and \ref{theoremtaboo} are essentially based 
on the assumption that $X(t)$ is a skip-free Markov process, and thus are analogous to  
similar results for other families of Markov processes such as simple random walks 
(see Mohanty and Panny \cite{MoPa}) and time-homogeneous diffusion processes 
(see Giorno~{\em et al.}\ \cite{GiNoRi}). 
\begin{example} \rm
Let $\{X(t),t\geq 0\}$ be the truncated BD process over $\{0,1,\ldots,N\}$, 
with $0$ and $N$ absorbing endpoints and rates 
$$
	\lambda_0=\lambda_N=\mu_0=\mu_N=0,
	\qquad\qquad \lambda_n=\lambda,
	\quad \mu_n=\mu
	\qquad(n=1,2,\ldots,N-1),
$$
with $\lambda,\mu>0$.  
For $k,n\in\{1,2,\ldots,N-1\}$ and $t>0$ the transition probabilities of $X(t)$ 
are given by (see formula (33) of B\"{o}hm and Mohanty~\cite{BoMo} 
amended of a misprint): 
\begin{equation}
	p_{k,n}(t)=e^{-(\lambda+\mu)t}\,\left(\frac{\lambda}{\mu}\right)^{(n-k)/2}
	\sum_{j=-\infty}^{+\infty}\left\{I_{n-k-2jN}-I_{n+k-2(j+1)N}\right\},
 \label{equation:19}
\end{equation}
where 
$$
	I_k:= I_k\left(2t\sqrt{\lambda\mu}\right)
	=\sum_{i=0}^{+\infty}\,\frac{\left(t\sqrt{\lambda\mu}\right)^{k+2i}}{ i!\,(k+i)!}
$$
denotes the modified Bessel function of first kind. Since $X(t)$ is symmetric, 
Proposition~\ref{theoremg} allows to obtain the FPT density from state $k$ 
to state $s$ when $N=2s$, with $k$ and $s$ integers: 
\begin{equation}
	g^+_{k,s}(t)=\frac{e^{-(\lambda+\mu)t}}{ t}\, 
	\left(\frac{\lambda}{\mu}\right)^{(s-k)/2}
	\sum_{j=-\infty}^{+\infty}\left\{(s-k-4sj)\,I_{s-k-4sj}-(s+k-4sj)\,I_{s+k-4sj}\right\},
 \label{eq:gpiues1}
\end{equation}
for $0<k<s$ and $t>0$.  Figure \ref{fig:1} shows some plots of density (\ref{eq:gpiues1}). Moreover, 
for $N=2s$ Proposition~\ref{theoremtaboo} leads to the $s$-avoiding transition probabilities  
\begin{equation}
 p_{k,n}^{\langle s\rangle}(t) = e^{-(\lambda+\mu)t}\,
	\left(\frac{\lambda}{\mu}\right)^{(n-k)/2}
	\sum_{j=-\infty}^{+\infty}\left\{I_{n-k-4sj}-I_{n+k-4sj}
	-I_{n+k-2s(2j+1)}+I_{n-k-2s(2j+1)}\right\}
 \label{eq:pknes1}
\end{equation}
for $0< n, k<s$ or $s< n, k<2s$, and $t>0$. See Figure \ref{fig:2} for some plots of probabilities 
(\ref{eq:pknes1}). 
\begin{figure}[t]  
\begin{center}
\epsfxsize=7.5cm
\epsfbox{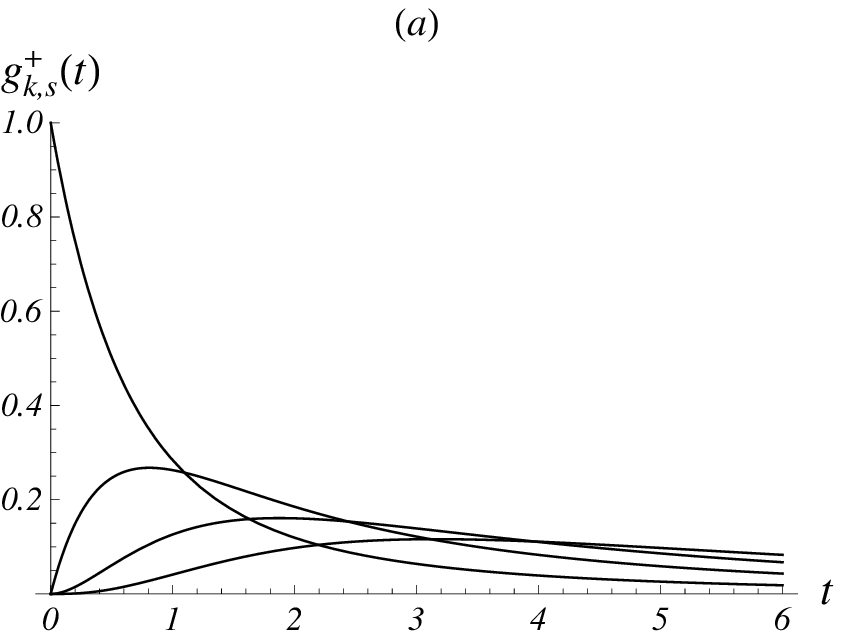}
$\;\;$
\epsfxsize=7.5cm
\epsfbox{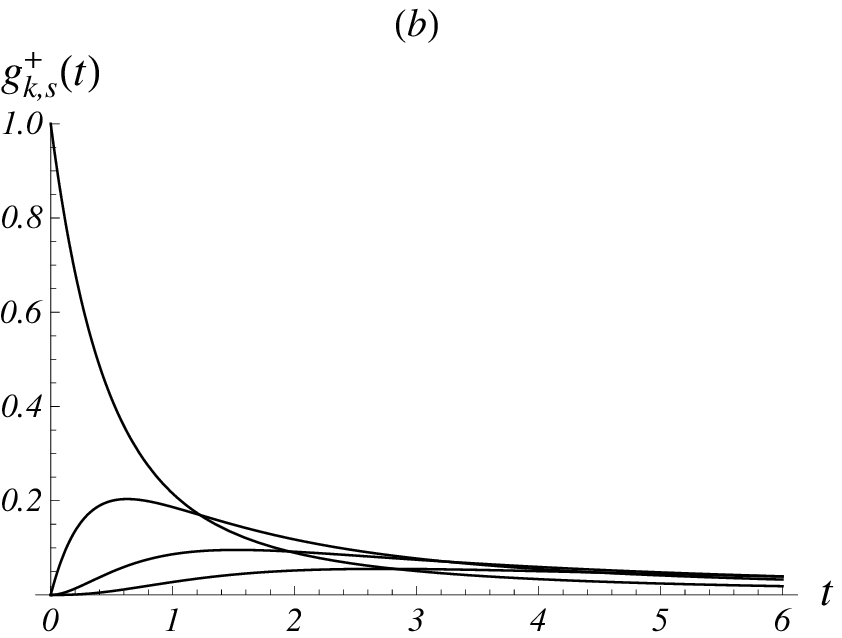}
\caption{The upward FPT density (\ref{eq:gpiues1}) for $\lambda=1$, $s=10$ and $k=6,7,8,9$ (from 
bottom to top near the origin),  with (a) $\mu=0.5$  and  (b) $\mu=1$.}
\label{fig:1}
\end{center}
\end{figure}
\begin{figure}[t]  
\begin{center}
\epsfxsize=7.5cm
\epsfbox{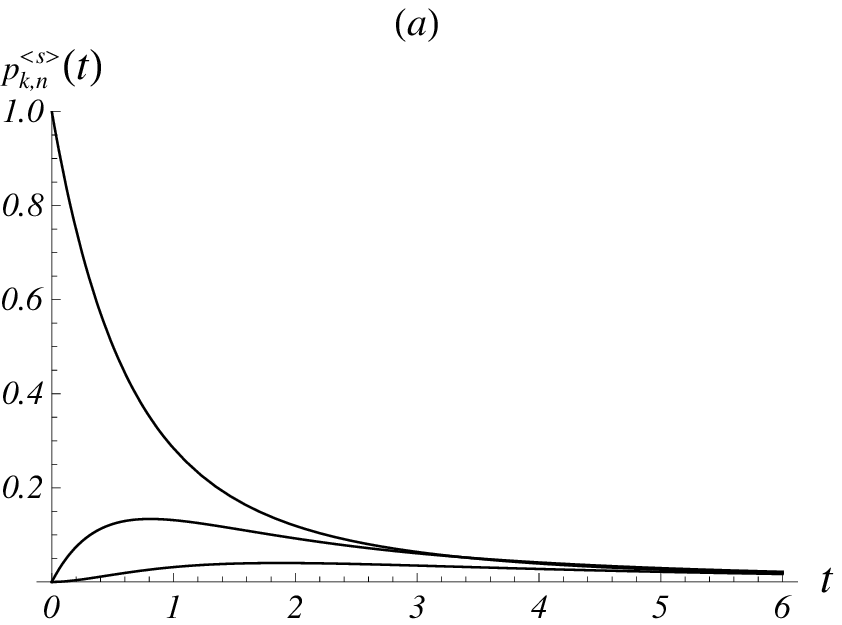}
$\;\;$
\epsfxsize=7.5cm
\epsfbox{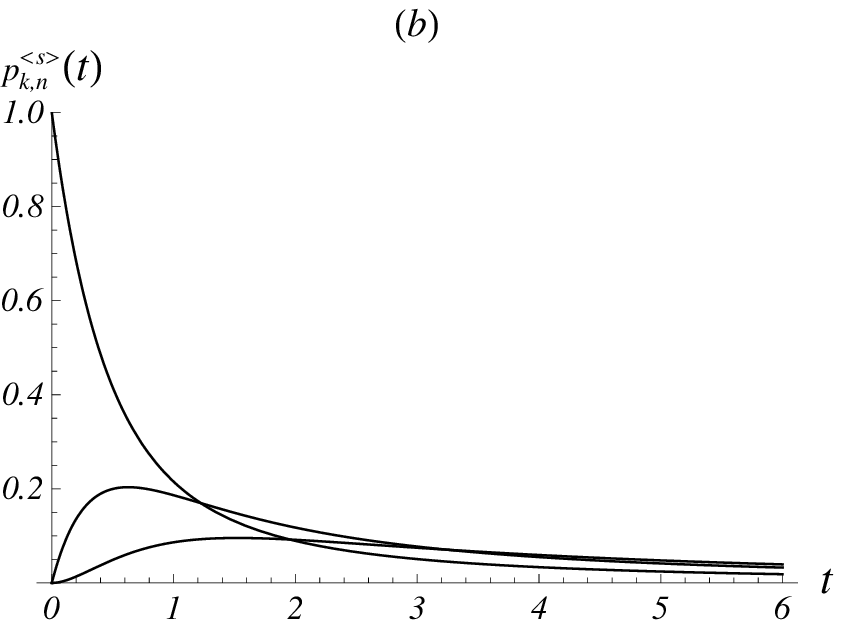}
\caption{The  $s$-avoiding transition probabilities (\ref{eq:pknes1}) for $\lambda=1$, $s=10$, 
$k=9$ and $n=7,8,9$ (from bottom to top near the origin),  with (a) $\mu=0.5$  and  (b) $\mu=1$.}
\label{fig:2}
\end{center}
\end{figure}
\hfill $\diamond$
\end{example}
%
\subsection{Reflecting endpoints}\label{section:4}
The symmetry properties treated in Section \ref{section:2} for BD processes with absorbing 
endpoints can be considered also in the case of reflecting endpoints, under similar hypotheses. 
\par
In this section we assume that  $\{X(t),t\geq 0\}$ is a BD process having state space $\cal S$, 
with $\mu_0=\lambda_N=0$, $\lambda_0,\mu_N>0$ and $\lambda_n,\mu_n>0$ for $n=1,2,\ldots,N-1$, 
so that the endpoints $0$ and $N>1$ are reflecting states. In this case, similarly to Theorem 
\ref{teorsimmassorb} one can prove that the transition probabilities of $X(t)$ satisfy the symmetry relation 
\begin{equation}	
	p_{k,n}(t)=p_{N-k,N-n}(t)\qquad (n=0,1,\ldots,N), 		
	\label{eq:symmrefl}
\end{equation}	
and we say that $X(t)$ is symmetric, if and only if 
$$
	\lambda_n=\mu_{N-n}\qquad (n=0,1,\ldots,N).
$$
In analogy to the case of absorbing endpoints, for symmetric BD processes with reflecting endpoints 
one can obtain various closed-form expressions. We limit ourselves to recall the following result, which is 
concerning the case $N=2s$ with $s$ integer (see Theorem 4.3 of \cite{DiCr98}). 
\begin{proposition}\label{teorgtaboorefl} 
Let the BD process $\{X(t),t\geq 0\}$ be   symmetric  in the sense of Eq.\ (\ref{eq:symmrefl}), with state 
space $\{0,1,$ $\ldots,2s\}$, where $0$ and $2s$ are reflecting states, with $s$ a positive integer. 
Then, one has:
\begin{eqnarray*}
	&& g^+_{k,s}(t)=\mu_{s+1}\,\big[p_{k,s-1}(t)-p_{k,s+1}(t)\big]
	\qquad\qquad (0\leq k<s), \\
 && p_{k,n}^{\langle s\rangle}(t) = p_{k,n}(t)-p_{2s-k,n}(t)	= p_{k,n}(t)-p_{k,2s-n}(t)
	\qquad (0\leq n, k<s\quad\hbox{or}\quad s<n, k\leq 2s). 
\end{eqnarray*}
\end{proposition}
\par
\begin{figure}[t]  
\begin{center}
\epsfxsize=7.5cm
\epsfbox{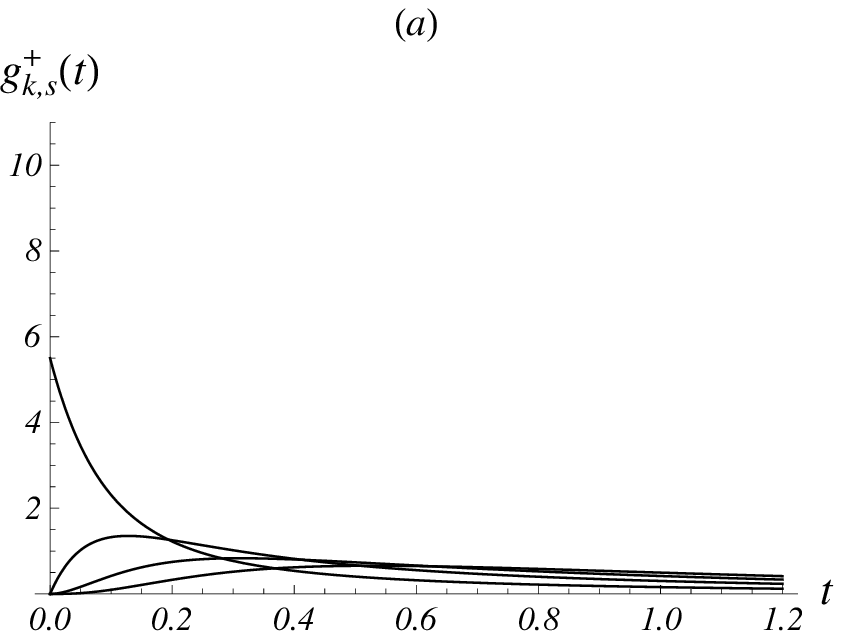}
$\;\;$
\epsfxsize=7.5cm
\epsfbox{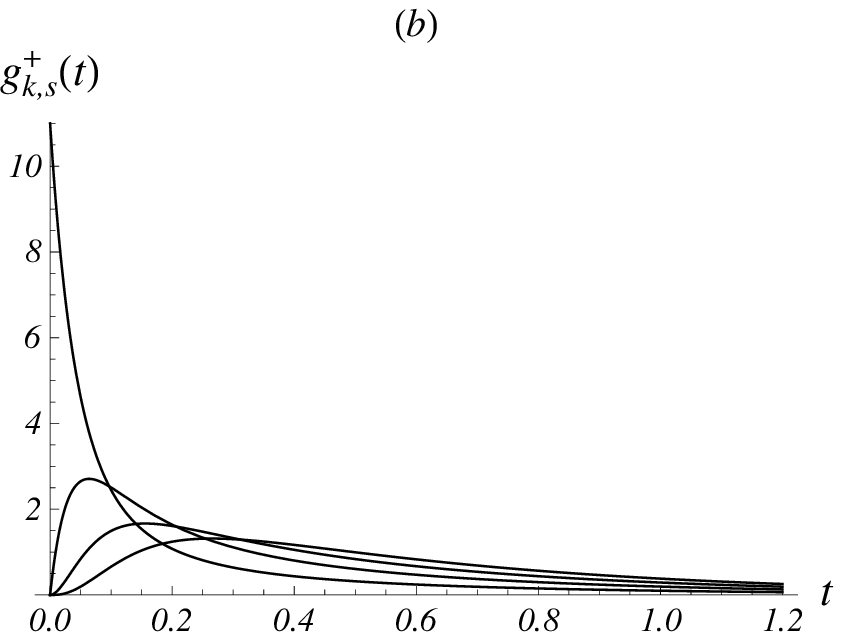}
\caption{The upward FPT density (\ref{eq:gcasorifl}) for  $s=10$ and $k=6,7,8,9$ (from 
bottom to top near the origin),  with (a) $\alpha =0.5$  and  (b) $\alpha =1$.}
\label{fig:3}
\end{center}
\end{figure}
\begin{figure}[t]  
\begin{center}
\epsfxsize=7.5cm
\epsfbox{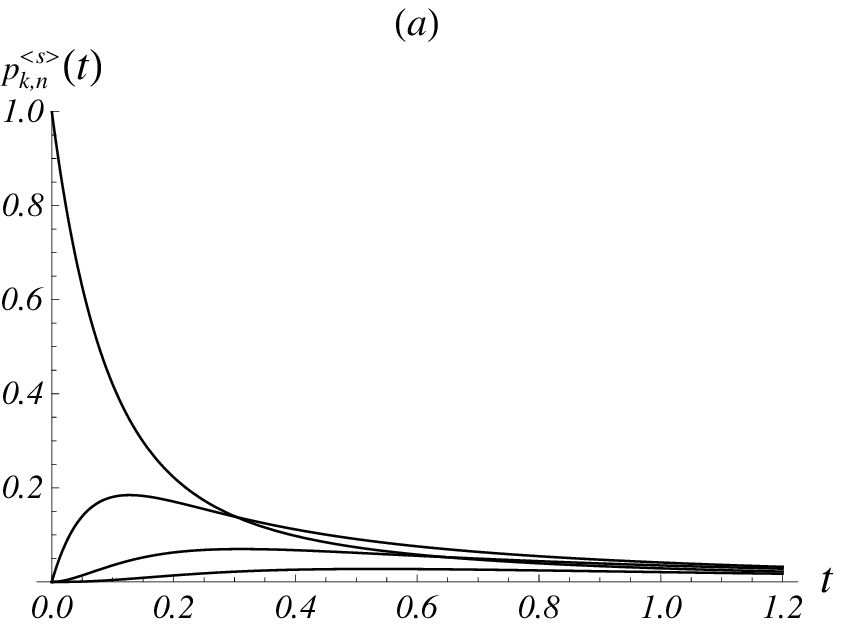}
$\;\;$
\epsfxsize=7.5cm
\epsfbox{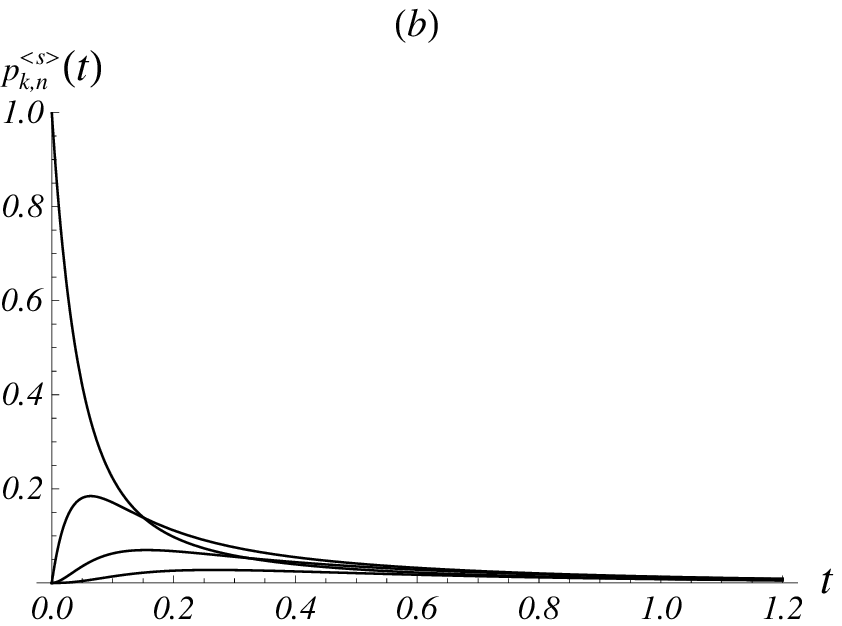}
\caption{The  $s$-avoiding transition probabilities (\ref{eq:pkncasorifl}) for  $s=10$, $k=9$ and 
$n=6,7,8,9$ (from bottom to top near the origin),  with (a) $\alpha=0.5$  and  (b) $\alpha=1$.}
\label{fig:4}
\end{center}
\end{figure}

Note that for small values of $t$, the FPT density $g^+_{k,s}(t)$ behaves as specified in Eq.\ (\ref{eq:gsmallt}). 
\begin{example}\rm
Let $\{X(t),t\geq 0\}$ be a  truncated BD process having 
state space $\{0,1,\ldots,N\}$ with reflecting endpoints and transition rates
$$
 \lambda_n=\alpha\,(N-n), \qquad 
 \mu_n=\alpha\,n \qquad (n=0,1,\ldots,N),
$$
where $\alpha>0$. This process is  symmetric  in the sense of Eq.\ (\ref{eq:symmrefl}), 
since its transition probabilities for $k,n\in\{0,1,\ldots,N\}$ are given by 
(see Giorno~{\em et al.}~\cite{GiNeNo}, for instance): 
\begin{equation}
	p_{k,n}(t)=\frac{1}{ 2^N}\sum_{j={\rm max}\{0,n+k-N\}}^{{\rm min}\{n,k\}}
	{k\choose j}{N-k\choose n-j}
	\big(1-e^{-2\alpha t}\big)^{n+k-2j}\,
	\big(1+e^{-2\alpha t}\big)^{N-(n+k-2j)}.
 \label{equation:20}
\end{equation}
Hence, if $N=2s$, with $s$ integer, from Proposition \ref{teorgtaboorefl} one has 
the following closed-form expressions of the FPT density through the symmetry state 
$s$: 
\begin{eqnarray}
	&& \hspace{-1cm} 
	g^+_{k,s}(t)=\frac{\alpha\,(s+1)}{ 2^{2s}}\,\sum_{j=0}^k{k\choose j}
	\Bigg[{2s-k\choose s-1-j}\big(1-e^{-2\alpha t}\big)^{s-1+k-2j}\,
	\big(1+e^{-2\alpha t}\big)^{s+1-k+2j}
	\nonumber \\
 && \;\; -{2s-k\choose s+1-j}\big(1-e^{-2\alpha t}\big)^{s+1+k-2j}\,
	\big(1+e^{-2\alpha t}\big)^{s-1-k+2j}\Bigg] 
	\qquad (0\leq k<s),
	\label{eq:gcasorifl}
\end{eqnarray}
and of the $s$-avoiding transition probabilities (for $0\leq n, k<s$ or $s<n, k\leq 2s$):
\begin{eqnarray}
 && p_{k,n}^{\langle s\rangle}(t) =\frac {1}{ 2^{2s}}\Bigg[
	\sum_{j={\rm max}\{0,n+k-2s\}}^{{\rm min}\{n,k\}}
	{k\choose j}{2s-k\choose n-j}
	\big(1-e^{-2\alpha t}\big)^{n+k-2j}\,
	\big(1+e^{-2\alpha t}\big)^{2s-n-k+2j}
	\nonumber \\
 && \qquad\qquad\;\; -\sum_{j={\rm max}\{0,n-k\}}^{{\rm min}\{n,2s-k\}}
	{2s-k\choose j}{k\choose n-j}
	\big(1-e^{-2\alpha t}\big)^{n+2s-k-2j}\,
	\big(1+e^{-2\alpha t}\big)^{k-n+2j}\Bigg]. 
	\label{eq:pkncasorifl}
\end{eqnarray}
Some plots of functions (\ref{eq:gcasorifl}) and (\ref{eq:pkncasorifl}) are shown in 
Figures \ref{fig:3} and \ref{fig:4}, respectively. 
\hfill $\diamond$
\end{example}
\par
Another example of symmetric  truncated BD process with reflecting endpoints is characterized by 
birth and death rates
$$
 \lambda_n=\alpha\,(N-n)^2, \qquad 
 \mu_n=\alpha\,n^2 \qquad (n=0,1,\ldots,N),
$$
with $\alpha>0$, and has been studied by Roehner and Valent \cite{RoVa}. 
\subsection{Bilateral processes}\label{section:bil}
It is worthwhile to point out that the symmetry properties of truncated BD processes 
with reflecting endpoints can be straightforwardly extended to bilateral BD processes. 
We recall that the term `bilateral' refers to stochastic processes defined on the whole set of 
integers (see, Pruitt \cite{Pr63}). Hence, we consider a bilateral BD process $\{X(t), t\geq 0\}$
having state space $\mathbb Z$, with birth and death rates $\lambda_n,\mu_n>0$ for all 
$n\in \mathbb Z$. Results similar to those given in Section \ref{section:4} can be stated in this 
case by noting that the transition probabilities of $X(t)$ satisfy the symmetry relation 
\begin{equation}	
	p_{k,n}(t)=p_{-k,-n}(t)\qquad (k,n\in \mathbb Z), 		
	\label{eq:symmbil}
\end{equation}	
for all $t>0$, if and only if 
$$
	\lambda_n=\mu_{-n}\qquad (n\in \mathbb Z).
$$
Note that here, for brevity, we consider the case in which 0 is the symmetry state. In analogy with 
Proposition \ref{teorgtaboorefl}, thanks to the symmetry property (\ref{eq:symmbil}) we have 
the following 
\begin{proposition}\label{teorbil} 
Let the bilateral BD process $\{X(t),t\geq 0\}$ be symmetric in the sense of Eq.\ (\ref{eq:symmbil}). 
Then, the upward and downward FPT densities through state 0, for $t>0$, are given by: 
\begin{align*}
	& g^+_{k,0}(t) =\mu_1 \left[p_{k,-1}(t)- p_{k,1}(t)\right]
	\qquad (k<0), 
	\\
	& g^-_{k,0}(t) =\mu_1 \left[p_{k,1}(t)  -  p_{k,-1}(t)\right]
	\qquad (k>0).
\end{align*}
Furthermore, the 0-avoiding transition probability of $X(t)$ can be expressed as
\begin{align*}
	p_{k,n}^{\langle 0\rangle}(t) = p_{k,n}(t)-p_{-k,n}(t)	 
	\qquad (n, k<0\quad\hbox{or}\quad 0<n, k),
\end{align*}
for all $t>0$. 
\end{proposition}
\par
We conclude this section by showing certain bilateral BD processes possessing the symmetry property 
considered in (\ref{eq:symmbil}).
\begin{description}
\item{\em (i)} 
First  we consider the BD process having sigmoidal-type rates 
\begin{equation*}
\lambda_n=\lambda\,\displaystyle\frac{1+c \Big(\displaystyle\frac{\mu}{\lambda}\Big)^{n+1}}
{1+c\Big(\displaystyle\frac{\mu}{\lambda}\Big)^{n}},  
 \qquad 
 {\mu_n}=\mu\,\displaystyle\frac{1+c \Big(\displaystyle\frac{\mu}{\lambda}\Big)^{n-1}}
{1+c \Big(\displaystyle\frac{\mu}{\lambda}\Big)^{n}} 
 \qquad  (n\in\mathbb{Z}),
\end{equation*}
with $\lambda,\mu>0$ and $c\geq 0$. Generally, the transition probabilities of this process exhibits 
a bimodality. See Hongler and Parthasarathy \cite{HoPa2008} and Di Crescenzo and Martinucci 
\cite{DiCrMa2009} for the symmetry property and other results. 
\item{\em (ii)} 
Other bilateral BD processes of interest are characterized by birth and death rates of alternating type. 
The process with rates 
\begin{equation}
\lambda_n=\left\{
\begin{array}{ll}
\lambda, & n \hbox{ even} \\
\mu, & n \hbox{ odd} 
\end{array}
\right.
 \qquad 
 {\mu_n}=\left\{
\begin{array}{ll}
\mu, & n \hbox{ even} \\
\lambda, & n \hbox{ odd} 
\end{array}
\right. 
 \qquad  (n\in\mathbb{Z}),
 \label{eq:alt1}
\end{equation}
for $\lambda,\mu>0$, has been studied by Conolly {\em et al.\/} \cite{CoPaDh} and 
Tarabia {\em et al.\/} \cite{TaTaEB09}. A suitable modification of the previous stochastic 
model yields the following birth and death rates:
\begin{equation}
\lambda_n=\left\{
\begin{array}{ll}
\lambda, & n \hbox{ even} \\
\mu, & n \hbox{ odd} 
\end{array}
\right.
 \qquad 
 {\mu_n}=\left\{
\begin{array}{ll}
\lambda, & n \hbox{ even} \\
\mu, & n \hbox{ odd} 
\end{array}
\right. 
 \qquad  (n\in\mathbb{Z}),
  \label{eq:alt2}
\end{equation}
with $\lambda,\mu>0$. Various symmetries and other properties of this process have been 
investigated in Di Crescenzo {\em et al.\/} \cite{DiCrIuMa12}, whereas a suitable extension 
of the models (\ref{eq:alt1}) and (\ref{eq:alt2}) can be found in \cite{DiCrMaMa14}.
\end{description}
%
\section{\bf Symmetry properties of bilateral processes with catastrophes}\label{section:3}
In this section we shall consider bilateral BD processes subject to catastrophes, i.e.\ jumps toward 
state $0$  occurring randomly in time. Even in this case the spatial symmetry can be used to obtain 
closed-form results related to the FPT problem through state 0. We notice that certain bilateral BD 
processes subject to catastrophes have been recently employed to describe double-ended 
queueing systems (cf.\ Di Crescenzo {\em et al.}\ \cite{DiCrGiKKNo12}). Various recent results 
on BD processes with catastrophes can be found in Dimou and Economou \cite{DimEco13}, 
Giorno and Nobile \cite{GiNo}, Giorno {\em et al.}\ \cite{GiNoSp} and Zeifman {\em et al.}\ \cite{ZeSaPa}. 
\par
Specifically, hereafter we refer to a BD process with jumps $\{X(t), t\geq 0\}$, having state space 
$\mathbb Z$. It is a continuous-time Markov chain characterized by the following transitions
concerning births, deaths and catastrophes, respectively: 
\par
{\em (a)} \ from $n\in\mathbb Z$ to $n+1$, with rate $\lambda_n$, 
\par
{\em (b)} \ from $n\in\mathbb Z$ to $n-1$, with rate $\mu_n$, 
\par
{\em (c)} \ from $n\in\mathbb{Z}\setminus \{0\}$ to $0$, with rate $\alpha_n$. 
\\
Hence, for all $t>0$ the transition probabilities $p_{k,n}(t)$ satisfy the following system:
\begin{align*}  
 & \frac{d}{dt} p_{k,n}(t)  
 = -(\lambda_n+\mu_n+\alpha_n)\,p_{k,n}(t)+\lambda_{n-1}\,p_{k,n-1}(t)+\mu_{n+1}\,p_{k,n+1}(t), 
 \qquad n\in\mathbb{Z} \setminus\{0\}, 
 \nonumber
 \\
 & \frac{d}{dt} p_{k,0}(t)
 = -(\lambda_0+\mu_0)\,p_{k,0}(t)+\lambda_{-1}\,p_{k,-1}(t)+\mu_1\,p_{k,1}(t)
 +\sum_{r\in\mathbb Z\setminus \{0\}} \alpha_r\,p_{k,r}(t).
\end{align*}
It is not hard to see that  $X(t)$ has a central 
symmetry with respect to state $0$, i.e.
\begin{equation}
 p_{-k,-n}(t)=p_{k,n}(t) \qquad \hbox{for all $t>0$ and $k,n\in \mathbb{Z}$,}
 \label{eq:simmbilatj}
\end{equation}
if and only if 
$$
 \lambda_n=\mu_{-n}\quad \hbox{for all $n\in \mathbb{Z}$} \qquad \hbox{and}
 \qquad \alpha_n=\alpha_{-n} \quad\hbox{for all $n\in \mathbb{Z}\setminus \{0\}$.} 
$$ 
In this case the process $X(t)$ is not skip-free due to the presence of catastrophes toward state 0, and 
thus its sample paths are no more `continuos'. Nevertheless, since all sample paths from a negative state 
to a positive one (and vice versa) are forced to pass through 0, the symmetry-based approach with 
respect to state 0 is still valid. In this case, a relevant role is played by the {\em probability currents\/} 
in state 0, defined as 
\begin{align*}
	 h^+_{k,0}(t) &=\lim_{\tau\downarrow 0} \frac{1}{\tau}P\{X(t+\tau)=0, X(t)<0\,|\,X(0)=k\} \\
	&=\lambda_{-1}\,p_{k,-1}(t)  +\sum_{n<0} \alpha_n \, p_{k,n}(t),
	\\
	 h^-_{k,0}(t) &=\lim_{\tau\downarrow 0} \frac{1}{\tau}P\{X(t+\tau)=0, X(t)>0\,|\,X(0)=k\} \\
	&=\mu_1\,p_{k,1}(t)  +\sum_{n>0} \alpha_n \, p_{k,n}(t),
\end{align*}
for $k\in  \mathbb{Z}$ and $t>0$. Similarly as Proposition \ref{teorgtaboorefl}, from (\ref{eq:simmbilatj}) 
one can obtain the following results. 
\begin{proposition}\label{teorcatast} 
Let $\{X(t),t\geq 0\}$ be a bilateral BD process with catastrophes. If it is symmetric  in the sense of 
Eq.\ (\ref{eq:simmbilatj}), then the upward and downward FPT densities through 
state 0, for $t>0$, are given by:
\begin{align*}
	g^+_{k,0}(t)& =h^+_{k,0}(t)-h^-_{k,0}(t)
	\\
	& =\mu_1 \left[p_{k,-1}(t)- p_{k,1}(t)\right]
	+\sum_{n<0} \alpha_n \, p_{k,n}(t)-\sum_{n>0} \alpha_n \, p_{k,n}(t)
	\qquad (k<0), 
	\\
	g^-_{k,0}(t)& =h^-_{k,0}(t)-h^+_{k,0}(t)
	\\
	& =\mu_1 \left[p_{k,1}(t)  -  p_{k,-1}(t)\right]
	+\sum_{n>0} \alpha_n \, p_{k,n}(t)-\sum_{n<0} \alpha_n \, p_{k,n}(t)
	\qquad (k>0),
\end{align*}
and satisfy the following symmetry relation:
$$
	g^+_{k,0}(t)=g^-_{-k,0}(t) 	\qquad (k<0).
$$
Moreover, the 0-avoiding transition probability of $X(t)$, for $t>0$, can be expressed as:
\begin{align*}
	p_{k,n}^{\langle 0\rangle}(t) = p_{k,n}(t)-p_{-k,n}(t)	 
	\qquad (n, k<0\quad\hbox{or}\quad 0<n, k). 
\end{align*}
\end{proposition}
\par
We refer to Theorems 3.1 and 3.2 of Di Crescenzo  and  Nastro \cite{DiCrNa2004} for related results. 
%
\begin{example}\rm 
Let $X(t)$ be the bilateral birth-death process with catastrophes characterized by constant rates
$$
 \lambda_n=\lambda, \qquad 
 \mu_n=\mu, \qquad 
 \alpha_n=\alpha.
$$
As specified in \cite{DiCrNa2004}, for 
all $k,n\in \mathbb{Z}$ and $t>0$ the transition probabilities can be expressed as 
\begin{align*}
 p_{k,n}(t)
 &= e^{-\alpha t}\,\widehat{p}_{k,n}(t)
 + \alpha \int_0^t e^{ -\alpha \tau} \,\widehat{p}_{0,n}(\tau)\,d\tau
 \\ 
 &={\left(\frac{\lambda}{\mu}\right)}^{\!\!\frac{n-k}{2}}
 I_{n-k}\left(2 \sqrt{\lambda \mu} \,t\right) \, e^{-(\lambda+\mu+\alpha)t}
 + \alpha {\left(\frac{\lambda}{\mu}\right)}^{\!\!\frac{n}{2}} \!
 \int_0^t e^{ -(\lambda+\mu+\alpha)\tau}   
 I_n\left( 2 \sqrt{\lambda \mu} \,\tau\right)d\tau, 
\end{align*}
where $I_n(x)$ denotes the modified Bessel function of the first kind, and where 
\begin{equation}
 \widehat{p}_{k,n}(t):={\left(\frac{\lambda}{\mu}\right)}^{\!\!\frac{n-k}{2}}
 I_{n-k}(2 \sqrt{\lambda \mu} \,t) \,e^{-(\lambda+\mu)t}
 \label{equation:37}
\end{equation}
is the transition probability of the bilateral Poisson birth-death process with 
birth rate $\lambda$ and death rate $\mu$. Process $X(t)$ is symmetric in the sense 
of Eq.\ (\ref{eq:simmbilatj}). Hence, making use of Proposition \ref{teorcatast} 
if $\lambda=\mu$, for all $t>0$ and $k=1,2,\ldots$ we have the FPT density: 
\begin{eqnarray*}
  g^-_{k,0}(t) =e^{- (2\lambda+\alpha) t} 
 \Bigg\{\lambda\,\left[I_{k-1}(2\lambda\,t)-I_{k+1}(2\lambda\,t)\right]  
 + \alpha \sum_{j=1}^{+\infty} \left[I_{k-j}(2\lambda\,t) - I_{k+j}(2\lambda\,t) \right]\Bigg\}. 
\end{eqnarray*}
Moreover,  the 0-avoiding transition probability is 
\begin{eqnarray*}
 p^{\langle 0\rangle}_{k,n}(t) 
 = e^{-(2\lambda+\alpha) t}\,\left[I_{n-k}(2\lambda\,t) -I_{n+k}(2\lambda\,t) \right],
\end{eqnarray*}
for $ t>0$ and $k,n=1,2,\ldots$. 
\hfill $\diamond$
\end{example}
\section{Symmetry properties of  two-dimensional processes}\label{section:2d}
In this section we exploit a symmetry-based approach for two-dimensional  BD processes with 
constant rates, by extending some of the results provided in the previous sections. 
We essentially refer to some contributions given in Di Crescenzo and Martinucci \cite{DiCrMa2008}. 
\par
Let $\{{\bf X}(t)=[X_1(t), X_2(t)];\,t\geq 0\}$ be a two-dimensional 
BD process with state space $\mathbb{Z}^2$, and transition probabilities  
\begin{equation}
P({\bf n}, t\,|\,{\bf k})
=P\{{\bf X}(t)={\bf n}\,|\,{\bf X}(0)={\bf k}\}, \qquad t\geq 0,
\label{equation:1}
\end{equation}
with  ${\bf n}=(n_1, n_2)\in\mathbb{Z}^2$ and ${\bf k}=(k_1, k_2)\in\mathbb{Z}^2$. 
Let us introduce the birth and death rates of ${\bf X}(t)$:
\begin{align*} 
\lambda_1 &= \lim_{s \downarrow 0}\frac{1}{s}
P\{X_1(t+s)=n_1+1,\,X_2(t+s)=n_2\,|\,{\bf X}(t)={\bf n}\},\\
\lambda_2 &= \lim_{s \downarrow 0}\frac{1}{s}
P\{X_1(t+s)=n_1,\,X_2(t+s)=n_2+1\,|\,{\bf X}(t)={\bf n}\},\\
\mu_1 &= \lim_{s \downarrow 0}\frac{1}{s}
P\{X_1(t+s)=n_1-1,\,X_2(t+s)=n_2\,|\,{\bf X}(t)={\bf n}\},\\
\mu_2 &= \lim_{s \downarrow 0}\frac{1}{s}
P\{X_1(t+s)=n_1,\,X_2(t+s)=n_2-1\,|\,{\bf X}(t)={\bf n}\}.
\end{align*}
Clearly, for $t>0$, the transition probabilities are solution of the following system:
\begin{align} 
\frac{{ d}}{{ d}t} \,P({\bf n}, t\,|\,{\bf k}) &= 
-(\lambda_1+\lambda_2+\mu_1+\mu_2)\,P({\bf n}, t\,|\,{\bf k}) 
\nonumber \\
&+ \lambda_1 \,P(n_1-1, n_2, t\,|\,{\bf k})
+\lambda_2 \,P(n_1, n_2-1, t\,|\,{\bf k}) 
\nonumber \\
&+ \mu_1 \,P(n_1+1, n_2, t\,|\,{\bf k})
+\mu_2 \,P(n_1, n_2+1, t\,|\,{\bf k}),  \qquad \forall {\bf n}\in\mathbb{Z}^2,
\label{equation:2}
\end{align}
with initial condition 
$$
P({\bf n}, 0\,|\,{\bf k})=\prod_{i=1}^{2} \delta_{n_i, k_i}.
$$
By making use of the probability generating function of ${\bf X}(t)$ one can prove that, for all $t>0$, 
the transition probabilities are:
\begin{equation}
P({\bf n}, t\,|\,{\bf k})=\prod_{i=1}^{2} {\rm e}^{-(\lambda_i+\mu_i)t} 
I_{n_i-k_i}(2 \sqrt{\lambda_i \mu_i} \,t) 
\Big(\frac{\lambda_i}{\mu_i}\Big)^{\frac{n_i-k_i}{2}}, 
 \qquad {\bf k},{\bf n}\in\mathbb{Z}^2.
\label{equation:5bis}
\end{equation}
From (\ref{equation:5bis}) the following quasi-symmetry property follows:  
\begin{proposition}\label{proposition:1}
If there exists a constant $\xi>0$ such that 
\begin{equation}
 \frac{\lambda_1}{\lambda_2}=\frac{\mu_2}{\mu_1}=\xi,
 \label{equation:10}
\end{equation}
then for all ${\bf n},{\bf k}\in\mathbb{Z}^2$ and all $r\in \mathbb{Z}$ we have
\begin{equation}
 P(n_2-r, n_1+r, t\,|\,k_2-r, k_1+r) = \xi^{n_2-k_2-n_1+k_1}\,P({\bf n}, t\,|\,{\bf k}). 
 \label{equation:3}
\end{equation}
\end{proposition}
\par
This result extends  the quasi-symmetry property given in Theorem \ref{teorsimmassorb} 
for one-dimensional truncated birth-death processes. In this case we deal with a spatial 
symmetry in the plane with respect to the straight line $x_2=x_1+r$. Namely, for each 
sample path of ${\bf X}(t)$ going from ${\bf k}$ to ${\bf n}$ there exists a symmetric path 
going from $(k_2-r, k_1+r)$ to $(n_2-r, n_1+r)$, where $r$ is a fixed integer. Eq.\ (\ref{equation:3}) 
thus expresses that the ratio of the probabilities of the two symmetric paths is time-independent. 
\par
This property is useful to obtain various results on the 
FPT problem of ${\bf X}(t)$ through straight-lines $x_2=x_1+r$. 
For a fixed $r\in\mathbb{Z}$, we denote by 
$$
T_{r}({\bf k})=\inf\{t\geq 0: X_2(t)=X_1(t)+r\}, 
\qquad {\bf X}(0)={\bf k}, \quad k_2\neq k_1+r,
$$
the FPT of ${\bf X}(t)$ through the straight line $x_2=x_1+r$, 
conditional on ${\bf X}(0)={\bf k}\in \mathbb{Z}^2$. Let 
$$
h_r(t\,|\,{\bf k}):=\frac{{ d}}{{ d}t} 
P\{T_r({\bf k})\leq t\,|\,{\bf X}(0)={\bf k}\}, \qquad t>0  
$$
be the corresponding probability density function. It is not hard to see that 
the following identity holds:
\begin{equation}
h_r(t\,|\,{\bf k})=\sum_{x\in\mathbb{Z}} g(x,x+r,t\,|\,{\bf k}),
\label{equation:12}
\end{equation}
where 
\begin{equation}
 g(x,x+r,t\,|\,{\bf k}):=\frac{\partial}{\partial t_1} 
 P\{T_r({\bf k})\leq t_1, {\bf X}(t_2)=(x,x+r)\,|\,{\bf X}(0)={\bf k}\} 
 \Big|_{t_1=t,t_2=t}
\label{equation:23}
\end{equation}
is the sub-density of the first-passage through line $x_2=x_1+r$ in state 
$(x,x+r)$ at time $t$, conditional on ${\bf X}(0)={\bf k}$. 
Clearly, for $n_2 \geq n_1+r$, $k_2<k_1+r$ (first-passage from below) 
and for $n_2 \leq n_1+r$, $k_2>k_1+r$ (first-passage from above) 
the following continuity equation holds:
\begin{equation}
P({\bf n}, t\,|\,{\bf k})=\int_{0}^{t}
\sum_{x\in\mathbb{Z}} g(x,x+r,\tau\,|\,{\bf k})\,
P({\bf n}, t-\tau\,|\,x,x+r)\,{ d}\tau.
\label{equation:4}
\end{equation}
For any fixed $r\in\mathbb{Z}$ we set 
\begin{equation}
 P^{\langle r\rangle}({\bf n}, t\,|\,{\bf k})
 =P\{{\bf X}(t)={\bf n}, T_r({\bf k})>t\,|\,{\bf X}(0)={\bf k}\}, 
 \qquad k_2\neq k_1+r,
 \label{equation:24}
\end{equation}
which expresses the probability of a sample path from ${\bf k}$ to ${\bf n}$ at time $t$ 
which does not touch the straight line $x_2=x_1+r$, conditional on ${\bf X}(0)={\bf k}$. 
By adopting a customary nomenclature in the field of Markov chains (see, for instance, 
Asmussen \cite{As03}) probability (\ref{equation:24}) is called `taboo probability'. 
\par
Thanks to the symmetry considered in Proposition \ref{proposition:1}, hereafter we can 
express the taboo probability in terms of transition probabilities. 
\begin{theorem}
Under the assumptions of Proposition \ref{proposition:1} we have 
\begin{equation}
 P^{\langle r\rangle}({\bf n}, t\,|\,{\bf k})
 =P({\bf n}, t\,|\,{\bf k})-{\xi}^{n_1+r-n_2}\,P(n_2-r,n_1+r,t\,|\,{\bf k}),
\label{equation:11}
\end{equation}
with $n_2< n_1+r$, $k_2<k_1+r$ or with $n_2> n_1+r$, $k_2>k_1+r$.
\label{teorema1}
\end{theorem}
\par
Let us now see that, as a further consequence of the symmetry property given in Proposition \ref{proposition:1}, 
the FPT densities can be suitably expressed in terms of transition probabilities.
\begin{theorem}\label{teorema2}
Under the assumptions of Proposition \ref{proposition:1}, for $k_2\neq k_1+r$ 
and $t>0$ we have 
\begin{equation}
 g(x,x+r,t\,|\,{\bf k})=\frac{|k_2-k_1-r|}{t} \,P(x,x+r,t\,|\,{\bf k}), 
 \qquad x\in \mathbb{Z},
 \label{equation:13}
\end{equation}
and 
\begin{equation}
 h_r(t\,|\,{\bf k})
 =\frac{|k_2-k_1-r|}{t} \, P\{X_2(t)=X_1(t)+r\,|\,{\bf X}(0)={\bf k}\}.
 \label{equation:6}
\end{equation}
\end{theorem}
\par
We remark that, due to Eqs.\ (\ref{equation:3}) and (\ref{equation:6}), under assumption (\ref{equation:10})  
the following quasi-symmetry relation holds:
\begin{equation}
 h_{r}(t\,|\,k_2-r,k_1+r)=\xi^{k_1+r-k_2} \, h_{r}(t\,|\,{\bf k}).
 \label{eq:proph}
\end{equation}
\par
Let us now introduce the quantity 
\begin{align}
 \pi_r({\bf k})
 & =P\{T_r({\bf k})<\infty \,|\,{\bf X}(0)={\bf k}\}
 \nonumber \\
 & =\int_{0}^{+\infty} h_{r}(t\,|\,{\bf k})\,{  d}t, 
 \qquad\qquad k_2\neq k_1+r,
 \label{equation:14}
\end{align}
which is the probability that process ${\bf X}(t)$ ultimately crosses the straight line 
$x_2=x_1+r$ starting from the initial value ${\bf X}(0)={\bf k}$. Due to relation (\ref{eq:proph}), 
under assumption (\ref{equation:10}) we have the following symmetry relation: 
$$
 \pi_r(k_2-r,k_1+r)=\xi^{k_1+r-k_2}\,\pi_r({\bf k}), 
 \qquad k_2\neq k_1+r.
$$
In conclusion we point out that, under the assumptions of Proposition \ref{proposition:1}, 
the first crossing probability  (\ref{equation:14}) is given by
$$
 \pi_r({\bf k})=\left\{
 \begin{array}{ll}
 \xi^{k_2-k_1-r} & \hbox{if }\lambda_1+\mu_2\geq \mu_1+\lambda_2, \; k_2<k_1+r,\\
 {}                     & \hbox{or }\lambda_1+\mu_2\leq \mu_1+\lambda_2, \; k_2>k_1+r,\\[0.2cm]
 1 & \hbox{otherwise.} 
 \end{array}
 \right.
$$

\subsection*{\bf Acknowledgments}
This paper is partially supported by GNCS-INdAM and Regione Campania (Legge 5).  

%
\end{document}